\documentclass[11pt]{article}

\usepackage{amsfonts,amssymb,amscd,amsmath,enumerate,verbatim,calc, times,url}
\usepackage{amsthm, psfrag,latexsym,epsfig,mdwlist,graphicx}

\theoremstyle{plain}
\newtheorem{theorem}{Theorem}[section]
\newtheorem{lemma}[theorem]{Lemma}
\newtheorem{proposition}[theorem]{Proposition}
\newtheorem{corollary}[theorem]{Corollary}
\theoremstyle{definition}
\newtheorem{definition}[theorem]{Definition}

\newtheorem{remark}[theorem]{Remark}
\newtheorem{example}[theorem]{Example}

\newcommand{\xs}{x_1,\ldots,x_n}            
\newcommand{\vs}{v_1,\ldots,v_n}            
\newcommand{\Ms}{M_1,\ldots,M_q}            
\newcommand{\Ns}{N_1,\ldots,N_s}            
\newcommand{\Fs}{F_1,\ldots,F_q}            
\newcommand{\dimn}{{\rm{dim}} \ }           
\newcommand{\F}{{\mathcal{F}}}              
\newcommand{\N}{{\mathcal{N}}}              
\newcommand{\I}{{\mathcal{I}}}              
\newcommand{\D}{\Delta}                     
\newcommand{\DN}{\Delta_N}                  
\newcommand{\al}{\alpha}                    
\newcommand{\be}{\beta}                     
\newcommand{\df}{\delta_{\mathcal{F}}}      
\newcommand{\dn}{\delta_{\mathcal{N}}}      
\newcommand{\lcm}{{\mathop{\rm{lcm}}}}      
\newcommand{\ddiv}{\ | \ }              
\newcommand{\ndiv}{\not | \ }              
\newcommand{\tuple}[1]{\langle #1 \rangle}  
\newcommand{\U}{{\mathcal{U}}}              
\newcommand{\PP}{{\mathcal{P}}}             
\newcommand{\rmv}[1]{\setminus \langle #1\rangle}
\newcommand{\x}[2]{x_{{#1},{#2}}}           
\newcommand{\height}{ {\rm{height}} \ }     
\newcommand{\alx}[1]{{#1{}^{\vee}}}         
\newcommand{\cmp}[1]{{#1}^c}                
\newcommand{\st}{\ | \ }
\newcommand{\erase}[1]{}                    
\newcommand{\vv}{\mbox{\boldmath{$\nu$}}}   
\newcommand{\nus}{\nu_1, \ldots, \nu_v}     
\newcommand{\QQ}{{\mathcal{Q}}}             
\newcommand{\RR}{{\mathcal{R}}}             
\newcommand{\pp}{{\mathfrak{p}}}             
\newcommand{\qq}{{\mathfrak{q}}}             
\newcommand{\ass}[1]{{\rm{Ass}}(#1)}         
\newcommand{\ann}[1]{{\rm{Ann}}(#1)}         
\newcommand{\scm}{sequentially Cohen-Macaulay }

\renewcommand{\wp}{\pp}
\newcommand{\QQs}[1]{\QQ_{\ #1}}

\date{2005}

\author{Sara Faridi\thanks{This research was
supported by NSERC, and by Laboratoire de combinatoire et
d'informatique math\'{e}matique at UQ\`AM, Montreal.}}

\title{\Large \sc Monomial ideals via square-free monomial ideals}

\begin{document}

\maketitle


\maketitle

\section{Introduction} 

In this paper we study monomial ideals using the operation
``polarization'' to first turn them into square-free monomial
ideals. Various forms of polarization appear throughout the literature
and have been used for different purposes in algebra and algebraic
combinatorics (for example, Weyman~\cite{W}, Fr\"oberg~\cite{Fr},
Schwartau~\cite{Schw}, or Rota and Stein~\cite{RS}). One of the most
useful features of polarization is that the chain of substitutions
that turn a given monomial ideal into a square-free one can be
described in terms of a regular sequence (Fr\"oberg~\cite{Fr}). This
fact allows many properties of a monomial ideal to transfer to its
polarization. Conversely, to study a given monomial ideal, one could
examine its polarization. The advantage of this latter approach is
that there are many combinatorial tools dealing with square-free
monomial ideals. One of these tools is Stanley-Reisner theory:
Schwartau's thesis \cite{Schw} and the book by St\"uckrad and Vogel
\cite{SV} discuss how the Stanley-Reisner theory of square-free
monomial ideals produces results about general monomial ideals using
polarization. Another tool for studying square-free monomial ideals,
which will be our focus here, is facet ideal theory, developed by the
author in \cite{F1}, \cite{F2} and \cite{F3}.

This paper is organized as follows. In Section~\ref{s:polarizaton} we
define polarization and introduce some of its basic properties. In
Section~\ref{s:facet} we introduce facet ideals and its features
that are relevant to this paper. In particular, we introduce simplicial
trees, which correspond to square-free monomial ideals with
exceptionally strong algebraic properties.
Section~\ref{s:monomial-polarized} extends the results of facet ideal
theory to general monomial ideals. Here we study a monomial ideal $I$
whose polarization is a tree, and show that many of the properties of
simplicial trees hold for such ideals. This includes
Cohen-Macaulayness of the Rees ring of $I$ (Corollary~\ref{c:Rees}),
$I$ being sequentially Cohen-Macaulay (Corollary~\ref{c:tree-scm}),
and several inductive tools for studying such ideals, such as
localization (see Section~\ref{ss:general-trees}).

Appendix~\ref{appendix} is an independent study of primary
decomposition in a sequentially Cohen-Macaulay module. We demonstrate
how in a sequentially Cohen-Macaulay module $M$, every submodule
appearing in the filtration of $M$ can be described in terms of the
primary decomposition of the $0$-submodule of $M$. This is used to
prove Proposition~\ref{p:duality}.

\bigskip 

\section{Polarization \hfill\break}\label{s:polarizaton}
 
\begin{definition}\label{d:polarization}  Let $R=k[\xs]$ be a polynomial
 ring over a field $k$. Suppose $M={x_1}^{a_1} \ldots {x_n}^{a_n}$ is a
monomial in $R$. Then we define the \emph{polarization} of $M$ to be
the square-free monomial
$$\PP(M)=\x{1}{1}\x{1}{2}\ldots\x{1}{a_1}\x{2}{1}\ldots\x{2}{a_2}\ldots\x{n}{1}
\ldots\x{n}{a_n}$$ in the polynomial ring $S=k[\x{i}{j} \ | \ 1 \leq i \leq
n, 1 \leq j \leq a_i]$.

If $I$ is an ideal of $R$ generated by monomials $\Ms$, then the
\emph{polarization} of $I$ is defined as: $$\PP(I)=\big (\PP(M_1),
\ldots, \PP(M_q)\big )$$ which is a square-free monomial ideal in a
polynomial ring $S$.
\end{definition}

 Here is an example of how  polarization works.

      \begin{example}\label{e:main-polar-example} Let $J=({x_1}^2,x_1x_2,{x_2}^3)
       \subseteq R=k[x_1,x_2]$. Then $$\PP(J)=
       (\x{1}{1}\x{1}{2},\x{1}{1}\x{2}{1},\x{2}{1}\x{2}{2}\x{2}{3})$$
       is the polarization of $J$ in the polynomial ring
       $$S=k[\x{1}{1},\x{1}{2},\x{2}{1},\x{2}{2},\x{2}{3}]$$
      \end{example}

Note that by identifying each $x_i$ with $\x{i}{1}$, one can consider
$S$ as a polynomial extension of $R$. Exactly how many variables $S$
has will always depend on what we polarize. Therefore, as long as we
are interested in the polarizations of finitely many monomials and
ideals, $S$ remains a finitely generated algebra.

Below we describe some basic properties of polarization, some of which
appear (without proof) in~\cite{SV}. Here we record the proofs where
appropriate.

\begin{proposition}[basic properties of polarization]\label{p:basic-polarization}Suppose that $R=k[\xs]$ is  a polynomial ring over a field $k$, and $I$
 and $J$ are two monomial ideals of $R$.

\begin{enumerate}

\item\label{i:1} $\PP(I+J)=\PP(I) + \PP(J)$;

\item\label{i:2} For two monomials $M$ and $N$ in $R$, $M \ddiv N$ if and only if
  $\PP(M) \ddiv \PP(N)$;

\item\label{i:3} $\PP(I \cap J)=\PP(I) \cap \PP(J)$;

\item\label{i:4} If $\pp=(x_{i_1}, \ldots, x_{i_r})$ is a (minimal)
prime containing $I$, then $\PP(\pp)$ is a (minimal) prime containing
$\PP(I)$; 

\item\label{i:5} If $\pp'=(\x{i_1}{e_1}, \ldots, \x{i_r}{e_r})$ is a
prime over $\PP(I)$, then $\pp=(x_{i_1}, \ldots, x_{i_r})$ is a prime
over $I$. Moreover, if $\pp'$ has minimal height (among all primes
containing $\PP(I)$), then $\pp$ must have minimal height as well (among
all primes containing $I$);

\item\label{i:6} $\height I = \height \PP(I)$;

\end{enumerate}

\end{proposition}

\begin{proof} \begin{enumerate}

\item Follows directly from Definition~\ref{d:polarization}.

\item Suppose that $M={x_1}^{b_1}\ldots {x_n}^{b_n}$ and $N={x_1}^{c_1}\ldots
{x_n}^{c_n}$, and suppose that 
$$\PP(M)=\x{1}{1}\ldots\x{1}{b_1}\ldots\x{n}{1} \ldots\x{n}{b_n}$$
 and $$\PP(N)=\x{1}{1}\ldots\x{1}{c_1}\ldots \x{n}{1}
\ldots\x{n}{c_n}.$$

If $M \ddiv N$, then $b_i \leq c_i$ for all $i$, which implies that
$\PP(M) \ddiv \PP (N)$. The converse is also clear using the same
argument.

\item Suppose that $I= (\Ms)$ and $J=(\Ns)$ where the generators are
           all monomials.  If $U={x_1}^{b_1}\ldots {x_n}^{b_n}$ is a
           monomial in $I\cap J$, then for some generator $M_i$ of $I$
           and $N_j$ and $J$, we have $M_i \ddiv U$ and $N_j \ddiv U$,
           hence by part ~\ref{i:2}, $\PP(M_i) \ddiv \PP(U)$ and
           $\PP (N_j) \ddiv \PP(U)$, which implies that $\PP(U) \in
           \PP(I) \cap \PP(J)$.
      
     Conversely, if $U'$ is a monomial in $\PP(I) \cap \PP(J)$, then
     for some generator $M_i={x_1}^{b_1}\ldots {x_n}^{b_n}$ of $I$ and
     $N_j={x_1}^{c_1}\ldots {x_n}^{c_n}$ and $J$ we have $\PP(M_i)
     \ddiv U'$ and $\PP (N_j) \ddiv U'$. This means that $\lcm
     (\PP(M_i), \PP(N_j)) \ddiv U'$. It is easy to see (by an argument
     similar to the one in part~\ref{i:2}) that $\lcm (\PP(M_i),
     \PP(N_j)) = \PP(\lcm(M_i, N_j))$. Since $\lcm(M_i, N_j)$ is one of
     the generators of $I \cap J$, it follows that $\PP(\lcm(M_i,
     N_j))$ is a generator of $\PP(I\cap J)$ and hence $U' \in
     \PP(I\cap J)$.
        
\item If $\pp=(x_{i_1}, \ldots, x_{i_r})$ is a minimal prime over
      $I=(\Ms)$, then for each of the $x_{i_j}$ there is a $M_t$ such
      that $x_{i_j} \ddiv M_t$, and no other generator of $\pp$ divides
      $M_t$. The same holds for the polarization of the two ideals:
      $\PP(\pp)=(\x{i_1}{1}, \ldots, \x{i_r}{1})$ and $\PP(I)=(\PP(M_1),
      \ldots, \PP(M_t))$, and so $\PP(\pp)$ is minimal over $\PP(I)$.

\item Suppose that $\pp'= (\x{i_1}{e_1},\ldots, \x{i_r}{e_r})$ is a
      prime lying over
      $\PP(I)$. Then for every generator $M_t$ of $I$, there is a
      $\x{i_j}{e_j}$ in $\pp'$ such that $\x{i_j}{e_j} \ddiv
      \PP(M_t)$. But this implies that $x_{i_j}\ddiv M_t$, and
      therefore $I \subseteq \pp=(x_{i_1},\ldots, x_{i_r})$.  

      Now suppose that $\pp'$ has minimal height $r$ over $\PP(I)$, and
      there is a prime ideal $\qq$ over $I$ with $\height \qq <r$. This
      implies (from part~\ref{i:4}) that $\PP(\qq)$, which is a prime of
      height less than $r$, contains $\PP(I)$, which is a
      contradiction.

     \item This follows from parts~\ref{i:4} and~\ref{i:5}.
      
      \end{enumerate}
\end{proof}

\begin{example} It is not true  that every minimal prime 
of $\PP(I)$ comes from a minimal prime of $I$. For example, let
$I=(x_1^2, x_1x_2^2)$. Then
$$\PP(I)=(\x{1}{1}\x{1}{2},\x{1}{1}\x{2}{1}\x{2}{2}).$$ The ideal
$(\x{1}{2},\x{2}{1})$ is a minimal prime over $\PP(I)$, but the
corresponding prime $(x_1,x_2)$ is not a minimal prime of $I$
(however, if we had taken any minimal prime of minimal height of
$\PP(I)$, e.g. $(\x{1}{1})$, then the corresponding prime over $I$
would have been minimal; this is part~\ref{i:5} above).
\end{example}

For a monomial ideal $I$ in a polynomial ring $R=k[\xs]$ as above,
there is a unique irredundant irreducible decomposition of the form
$$I=\qq_1\cap \ldots \cap \qq_m$$ where each $\qq_i$ is a primary ideal
generated by powers of the variables $\xs$ (see~\cite[Theorem~5.1.17]{Vi2}).

\begin{proposition}[polarization and primary 
decomposition]\label{p:polar-primary}
Let $I$ be a monomial ideal in a polynomial ring $R=k[\xs]$, and let
$\PP(I)$ be the polarization of $I$ in $S=k[\x{i}{j}]$ as described in
{\rm Definition~\ref{d:polarization}}. 

\begin{enumerate} 

\item\label{i:7} If $I=({x_{i_1}}^{a_1}, \ldots, {x_{i_r}}^{a_r})$ where the
$a_j$ are positive integers, then
\[ \PP(I)=\bigcap_{\stackrel{\mbox{$\scriptstyle 1 \leq c_j \leq
a_j$}}{ \mbox{$\scriptstyle 1 \leq j \leq r$}}} (\x{i_1}{c_1}, \ldots,
\x{i_r}{c_r})\]

\item\label{i:8} If $I=(x_{i_1}, \ldots, x_{i_r})^m$, where $1 \leq
           i_1, \ldots, i_r \leq n$ and $m$ is a positive integer,
           then $\PP(I)$ has the following irredundant irreducible
           primary decomposition:
\[ \PP(I)= \bigcap_{\stackrel{\mbox{$\scriptstyle 1 \leq c_j \leq m$}}{
\mbox{$\scriptstyle \Sigma c_j\leq m+r-1$}}} (\x{i_1}{c_1}, \ldots,
           \x{i_r}{c_r})\]

\item\label{i:9} Suppose that $I=\qq_1\cap \ldots \cap \qq_m$ is the
unique irredundant irreducible primary decomposition of $I$, such that
for each $i=1,\ldots, m$, $$\qq_i=({x_1}^{a_1^i},\ldots,
{x_n}^{a_n^i}),$$ where the $a_j^i$ are nonnegative integers, and if
$a_j^i=0$ we assume that ${x_j}^{a_j^i}=0.$ 

Then $\PP(I)$ has the following irreducible primary decomposition
(some primes might be repeated).
\[ \PP(I)= \bigcap_{\mbox{$\scriptstyle 1 \leq i \leq m$}} 
\bigcap_{\stackrel{\mbox{$\scriptstyle 1 \leq c_j \leq a_j^i$}}{
\mbox{$\scriptstyle 1 \leq j \leq n$}}} (\x{1}{c_1}, \ldots,
\x{n}{c_n})\] where when $a_j^i=0$, we assume that $c_j=\x{j}{0}=0$.

\end{enumerate}

\end{proposition}

\begin{proof}

\begin{enumerate}

\item We know that $$\PP(I)=(\x{i_1}{1}\ldots \x{i_1}{a_1}, \ldots,
\x{i_1}{1}\ldots \x{i_1}{a_r}).$$ Clearly the minimal primes of
$\PP(I)$ are $(\x{i_1}{c_1}, \ldots, \x{i_r}{c_r})$ for all $c_j \leq
a_j$. This settles the claim.

\item Assume, without loss of generality, that $I=(x_1, \ldots,
      x_r)^m$. So we can write
      $$I=\big ( {x_1}^{b_1}\ldots {x_r}^{b_r}\ | \ 0 \leq b_i \leq
      m,\ b_1 + \cdots + b_r =m \big )$$ so that $$\PP(I)=\big (
      \x{1}{1}\ldots \x{1}{b_1}\ldots \x{r}{1}\ldots \x{r}{b_r}\ | \ 0
      \leq b_i \leq m,\ b_1 + \cdots + b_r =m \big ).$$

      We first show that $\PP(I)$ is contained in the intersection of
      the ideals of the form $(\x{1}{c_1}, \ldots, \x{r}{c_r})$ described
      above. It is enough to show this for each generator of
      $\PP(I)$. So we show that $$\U = \x{1}{1}\ldots \x{1}{b_1}\ldots
      \x{r}{1}\ldots \x{r}{b_r} \in \I= (\x{1}{c_1}, \ldots,
      \x{r}{c_r})$$ where $ 0 \leq b_i \leq m$, $b_1 + \cdots + b_r
      =m$, $1 \leq c_j \leq m$ and $c_1 + \cdots + c_r \leq m+r-1$.
      
      If for any $i$, $b_i \geq c_i$, then it would be clear that $\U
      \in \I$.

      Assume $b_i \leq c_i-1$ for $i=1,\ldots,r-1$.  It follows that
       
      $$ \begin{array}{ll} m-b_r & = b_1+\cdots+b_{r-1}\\ & \leq
         c_1+\cdots+ c_{r-1} -(r-1)\\ & \leq m + r -1 -c_r - (r-1)\\ &
         = m - c_r
       \end{array}$$  

      which implies that $b_r \geq c_r$, hence $\U \in \I$.

      So far we have shown one direction of the inclusion.
      
      To show the opposite direction, take any monomial $$\U \in
      \bigcap  (\x{1}{c_1}, \ldots, \x{r}{c_r})$$ where
      $1 \leq c_j \leq m$ and $c_1 + \cdots + c_r \leq m+r-1$.

      Notice that for some $i \leq r$, $\x{i}{1} \ddiv \U$; this is
      because $\U \in ( \x{1}{1}, \ldots, \x{r}{1})$.

      We write $\U$ as
      $$\U=\x{1}{1}\ldots \x{1}{b_1}\ldots \x{r}{1}\ldots
      \x{r}{b_r}\U'$$ where $\U'$ is a monomial, and the $b_i$ are
      nonnegative integers such that for each $j < b_i$, $\x{i}{j}
      \ddiv \U$ (if $\x{i}{1} \ndiv \U$ then set $b_i=0$). We need to
      show that it is possible to find such $b_i$ so that $b_1 +
      \cdots + b_r = m$.

      Suppose $b_1 + \cdots + b_r \leq m-1$, and $\x{i}{b_i+1} \ndiv
      \U$ for $1 \leq i \leq r$. Then $$b_1 + \cdots + b_r +r \leq
      m+r-1,$$ hence $$\U \in (\x{1}{b_1+1}, \ldots, \x{r}{b_r+1})$$
      implying that $\x{i}{b_i+1} \ddiv \U$ for some $i$, which is a
      contradiction.

      Therefore $b_1, \ldots, b_r$ can be picked so that they add up
      to $m$, and hence $\U \in \PP(I)$; this settles the opposite
      inclusion.

\item This follows from part~\ref{i:7} and
Proposition~\ref{p:basic-polarization} part~\ref{i:3}.

\end{enumerate}

\end{proof}

\begin{corollary}[polarization and associated primes]\label{c:polar-ass}
Let $I$ be a monomial ideal in a polynomial ring $R=k[\xs]$, and let
$\PP(I)$ be its polarization in $S=k[\x{i}{j}]$ as described in
{\rm Definition~\ref{d:polarization}}. Then $(x_{i_1}, \ldots, x_{i_r}) \in
{\rm Ass}_R(R/I)$ if and only if $(\x{i_1}{c_1}, \ldots, \x{i_r}{c_r})
\in {\rm Ass}_S(S/\PP(I))$ for some positive integers $c_1, \ldots,
c_r$.
\end{corollary}

\erase{
(February 4, 2017) Took out, from statement of corollary: ``
 Moreover, if $(\x{i_1}{c_1}, \ldots, \x{i_r}{c_r}) \in {\rm
Ass}_S(S/\PP(I))$, then $(\x{i_1}{b_1}, \ldots, \x{i_r}{b_r}) \in {\rm
Ass}_S(S/\PP(I))$ for all $b_j$ such that $1 \leq b_j \leq c_j$.''

 because, for example, after polarization the monomial ideal
 $(x^2,y^2) \cap (x,y^3,z)$ becomes $$(x_1,y_1) \cap (x_1,y_2) \cap
 (x_2,y_1) \cap (x_2,y_2) \cap (x_1,y_1,z_1)\cap (x_1,y_2,z_1)\cap
 (x_1,y_3,z_1)$$ which is equal to
 $$(x_1,y_1) \cap (x_1,y_2) \cap (x_2,y_1) \cap (x_2,y_2) \cap
 (x_1,y_3,z_1).$$ Now $ (x_1,y_3,z_1)$ is an associated prime of the
 polarized ideal, but $ (x_1,y_2,z_1)$ is not.
 }
\begin{example} 
       Consider the primary decomposition of $J=({x_1}^2,
       {x_2}^3,x_1x_2)$: $$J=(x_1, {x_2}^3) \cap ({x_1}^2, x_2).$$ By
       Proposition~\ref{p:polar-primary}, $\PP(J)=(\x{1}{1}\x{1}{2},
       \x{2}{1}\x{2}{2}\x{2}{3},\x{1}{1}\x{2}{1})$ will have primary
       decomposition
       $$\PP(J)=  (\x{1}{1},\x{2}{1}) \cap 
                  (\x{1}{1},\x{2}{2}) \cap 
                  (\x{1}{1},\x{2}{3}) \cap 
                  (\x{1}{2},\x{2}{1}). $$

\end{example}

A very useful property of polarization is that the final polarized
ideal is related to the original ideal via a regular sequence. The
proposition below, which looks slightly different here than the
original statement in~\cite{Fr}, states this fact.

\begin{proposition}[Fr\"oberg~\cite{Fr}]\label{p:polarization} Let $k$
be a field and $$R=k[\xs]/(\Ms),$$ where $\Ms$ are monomials in the
 variables $\xs$, and let $$N_1=\PP(M_1), \ldots, N_q=\PP(M_q)$$ be a
 set of square-free monomials in the polynomial ring $$S=k[\x{i}{j} \
 | \ 1 \leq i \leq n, 1 \leq j \leq a_i]$$ such that for each $i$, the
 variable $\x{i}{a_i}$ appears in at least one of the monomials $N_1,
 \ldots, N_q$. Then the sequence of elements
  \begin{eqnarray}\label{e:sequence}\x{i}{1} - \x{i}{j} {\rm \ where \ }
  1 \leq i \leq n {\rm \ and\ } 1 < j \leq a_i \end{eqnarray} forms a
  regular sequence in the quotient ring $$R'=S/(N_1,\ldots, N_q)$$ and
  if $J$ is the ideal of $R'$ generated by the elements in {\rm
  (\ref{e:sequence})}, then $$R=R'/J.$$ Moreover, $R$ is
  Cohen-Macaulay (Gorenstein) if and only if $R'$ is.
\end{proposition}

\begin{example} Let $J$ and $R$ be as in
 Example~\ref{e:main-polar-example}. According to
 Proposition~\ref{p:polarization}, the sequence
      $$\x{1}{1}-\x{1}{2},\ \x{2}{1}-\x{2}{2},\ \x{2}{1}- \x{2}{3}$$
      is a regular sequence in $S/\PP(J)$, and $$R/J = S/\big
      (\PP(J)+(\x{1}{1}-\x{1}{2}, \x{2}{1}-\x{2}{2},\x{2}{1}-
      \x{2}{3})\big ).$$ \end{example}

\bigskip

\section{Square-free monomial ideals as facet 
ideals \hfill\break}\label{s:facet} 

Now that we have introduced polarization as a method of transforming a
monomial ideal into a square-free one, we can focus on square-free
monomial ideals. In particular, here we are interested in properties
of square-free monomial ideals that come as a result of them being
considered as facet ideals of simplicial complexes. Below we review
the basic definitions and notations in facet ideal theory, as well as
some of the basic concepts of Stanley-Reisner theory. We refer the
reader to~\cite{BH}, \cite{F1}, \cite{F2}, \cite{F3}, and \cite{S} for
more details and proofs in each of these topics.

\begin{definition}[simplicial complex, facet, subcollection and more] 
   A \emph{simplicial complex} $\D$ over a set of vertices $V=\{
  \vs \}$ is a collection of subsets of $V$, with the property that
  $\{ v_i \} \in \D$ for all $i$, and if $F \in \D$ then all
  subsets of $F$ are also in $\D$ (including the empty set). An
  element of $\D$ is called a \emph{face} of $\D$, and the
  \emph{dimension} of a face $F$ of $\D$ is defined as $|F| -1$,
  where $|F|$ is the number of vertices of $F$.  The faces of
  dimensions 0 and 1 are called \emph{vertices} and \emph{edges},
  respectively, and $\dimn \emptyset =-1$.  The maximal faces of
  $\D$ under inclusion are called \emph{facets} of $\D$. The
  dimension of the simplicial complex $\D$ is the maximal
  dimension of its facets.

We denote the simplicial complex $\D$ with facets $\Fs$ by
$$\D = \tuple{\Fs}$$ and we call $\{ \Fs \}$ the \emph{facet set} of
$\D$. A simplicial complex with only one facet is called a
\emph{simplex}. By a \emph{subcollection} of $\D$ we mean a simplicial
complex whose facet set is a subset of the facet set of $\D$.
\end{definition}

\begin{definition}[connected simplicial complex] A simplicial 
  complex
  $\D= \langle F_1, \ldots,$ $ F_q \rangle$ is \emph{connected} 
  if for every pair $i,j$, $1
  \leq i < j \leq q$, there exists a sequence of facets
  $F_{t_1},\ldots,F_{t_r}$ of $\D$ such that $F_{t_1}=F_i$,
  $F_{t_r}=F_j$ and $F_{t_s} \cap F_{t_{s+1}} \neq \emptyset$ for
  $s=1,\ldots,r-1$.
\end{definition}

\begin{definition}[facet/non-face ideals and complexes]\label{d:facet-ideal} 
 Consider a polynomial ring $R=k[\xs]$ over a field $k$ and a set
of indeterminates $\xs$. Let $I=(\Ms)$ be an ideal in $R$, where $\Ms$
are square-free monomials that form a minimal set of generators for
$I$.

\begin{itemize}
 
\item The \emph{facet complex} of $I$, denoted by $\df(I)$, is the
 simplicial complex over a set of vertices $\vs$ with facets $\Fs$,
 where for each $i$, $F_i=\{v_j \st x_j|M_i, $ $\ 1 \leq j \leq n \}$.
 The \emph{non-face complex} or the \emph{Stanley-Reisner complex} of
 $I$, denoted by $\dn(I)$ is be the simplicial complex over a set of
 vertices $\vs$, where $\{v_{i_1},\ldots, v_{i_s}\}$ is a face of
 $\dn(I)$ if and only if $x_{i_1}\ldots x_{i_s} \notin I$.

\item Conversely, if $\D$ is a simplicial complex over $n$ vertices
 labeled $\vs$, we define the \emph{facet ideal} of $\D$, denoted by
 $\F(\D)$, to be the ideal of $R$ generated by square-free monomials
 $x_{i_1}\ldots x_{i_s}$, where $\{v_{i_1},\ldots, v_{i_s}\}$ is a
 facet of $\D$.  The \emph{non-face ideal} or the
 \emph{Stanley-Reisner ideal} of $\Delta$, denoted by $\N(\Delta)$, is
 the ideal of $R$ generated by square-free monomials $x_{i_1}\ldots
 x_{i_s}$, where $\{v_{i_1},\ldots, v_{i_s}\}$ is not a face of
 $\Delta$.

\end{itemize}

\end{definition}

Throughout this paper we often use a letter $x$ to denote both a
vertex of $\D$ and the corresponding variable appearing in $\F(\D)$,
and $x_{i_1}\ldots x_{i_r}$ to denote a facet of $\D$ as well as a
monomial generator of $\F(\D)$.

\begin{example}\label{main-example} If $\D$ is the simplicial complex
  $\tuple{xyz,yu,uvw}$ drawn below,
\[ \includegraphics{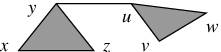} \]
  then $\F(\D)=(xyz,yu,uvw)$ and $\N(\D)=(xu,xv,xw, yv,
  yw,zu,zv,zw)$ are its facet ideal and nonface (Stanley-Reisner)
  ideal, respectively.
\end{example}

Facet ideals give a one-to-one correspondence between simplicial
complexes and square-free monomial ideals. 

Next we define the notion of a vertex cover. The combinatorial idea
here comes from graph theory. In algebra, it corresponds to prime
ideals lying over the facet ideal of a given simplicial complex.

\begin{definition}[vertex covering, independence, unmixed]\label{d:unmixed} 
  Let $\D$ be a simplicial complex with vertex set $V$. A \emph{vertex
  cover} for $\D$ is a subset $A$ of $V$ that intersects every facet of $\D$.
  If $A$ is a minimal element (under inclusion) of the set of vertex covers of
  $\D$, it is called a \emph{minimal vertex cover}. The smallest of the
  cardinalities of the vertex covers of $\D$ is called the \emph{vertex
  covering number} of $\D$ and is denoted by $\al(\D)$.  A simplicial complex
  $\D$ is \emph{unmixed} if all of its minimal vertex covers have the same
  cardinality.

  A set $\{F_1, \ldots, F_u\}$ of facets of $\D$ is called an
  \emph{independent set} if $F_i \cap F_j =\emptyset$ whenever $i \neq
  j$. The maximum possible cardinality of an independent set of facets
  in $\D$, denoted by $\beta(\D)$, is called the \emph{independence
  number} of $\D$. An independent set of facets which is not a proper
  subset of any other independent set is called a \emph{maximal
  independent set} of facets.
\end{definition}

\begin{example}\label{highlight-example}  If $\D$ is the simplicial 
complex in Example~\ref{main-example}, then the vertex covers of $\D$
  are:
$${\bf
\{x,u\},\{y,u\},\{y,v\},\{y,w\},\{z,u\}},\{x,y,u\},\{x,z,u\},\{x,y,v\},
\ldots .$$ The first five vertex covers above (highlighted in bold),
are the minimal vertex covers of $\D$. It follows that $\al(\D)=2$,
and $\D$ is unmixed. On the other hand, $\{xyz,uvw\}$ is the largest
maximal independent set of facets that $\D$ contains, and so
$\beta(\D)=2$.
\end{example}

\begin{definition}[Alexander dual]\label{d:alex}  Let $I$ be a square-free
 monomial ideal in the polynomial ring $k[V]$ with $V=\{\xs\}$, and
  let $\DN$ be the non-face complex of $I$ (i.e. $\DN=\dn(I)$). Then the
  \emph{Alexander dual} of $\DN$ is the simplicial complex
  $$\alx{\DN} = \{ F \subset V \st \cmp{F} \notin \DN\}$$ where
$\cmp{F}$ is the complement of the face $F$ in $V$.  

We call the nonface ideal of $\alx{\DN}$ the \emph{Alexander dual} of
$I$ and denote it by $\alx{I}$.\end{definition}

\subsection{Simplicial Trees}

Considering simplicial complexes as higher dimensional graphs, one can define
the notion of a \emph{tree} by extending the same concept from graph theory.
Before we define a tree, we determine what ``removing a facet'' from a
simplicial complex means. We define this idea so that it corresponds to
dropping a generator from the facet ideal of the complex.

\begin{definition}[facet removal]\label{d:removal}  Suppose $\D$
  is a simplicial complex with facets $F_1,$ $\ldots, F_q$ 
  and $\F(\D)=(\Ms)$ its
  facet ideal in $R=k[\xs]$. The simplicial complex obtained by
  \emph{removing the facet} $F_i$ from $\D$ is the simplicial
  complex
 $$\D \rmv{F_i}=\tuple{F_1,\ldots,\hat{F}_{i},\ldots,F_q}$$ and $\F(\D
\rmv{F_i}) = (M_1,\ldots,\hat{M}_{i} ,\ldots,M_q)$.
 \end{definition}
 
The definition that we give below for a simplicial tree is one generalized from
graph theory. See \cite{F1} and \cite{F2} for more on this concept.

\begin{definition}[leaf, joint] A facet $F$ of a simplicial complex 
is called a \emph{leaf} if either $F$ is the only facet of $\D$, or
  for some facet $G \in \D \rmv{F}$ we have
  $$F \cap (\D \rmv{F}) \subseteq G.$$

If $F \cap G \neq \emptyset$, the facet $G$ above is called a
\emph{joint} of the leaf $F$.
\end{definition}

Equivalently, a facet $F$ is a leaf of $\D$ if $F \cap (\D \rmv{F})$
is a face of $\D \rmv{F}$.

\begin{example}\label{leaf-example} {\rm Let $I=(xyz,yzu,zuv)$. Then $F=xyz$ is
 a leaf, but $H=yzu$ is not, as one can see in the picture below.

\[ \begin{tabular}{cccccc}
  $F \cap (\D \rmv{F}) =$& \includegraphics{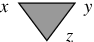} &$\cap$ &\includegraphics{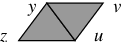} &= &
\includegraphics{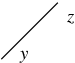} \\
&&&&\\
$H \cap (\D \rmv{H})=$& \includegraphics{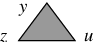} &$\cap$ &\includegraphics{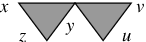} &= &
\includegraphics{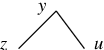}
\end{tabular} \]}
\end{example}

\begin{definition}[tree, forest]\label{d:tree}A connected simplicial 
  complex $\D$ is a \emph{tree} if every nonempty subcollection of
  $\D$ has a leaf.  If $\D$ is not necessarily connected, but
  every subcollection has a leaf, then $\D$ is called a forest.
\end{definition}

\begin{example}\label{free-example} The simplicial complexes in 
examples~\ref{main-example} and~\ref{leaf-example} are both trees, but
  the one below is not because it has no leaves. It is an easy
  exercise to see that a leaf must contain a free vertex, where a
  vertex is \emph{free} if it belongs to only one facet.
\[ \includegraphics{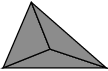} \]
\end{example}

One of the most powerful properties of simplicial trees from the point of view
of algebra is that they behave well under localization. This property makes it
easy to use induction on the number of vertices of a tree for proving its
various properties.

\begin{lemma}[Localization of a tree is a forest]\label{l:localization}
  Let $I \subseteq k[\xs]$ be the facet ideal of a simplicial tree,
  where $k$ is a field. Then for any prime ideal $\pp$ of $k[\xs]$, $\df
  (I_{\pp})$ is a forest. \end{lemma}

\begin{proof} See~\cite[Lemma 4.5]{F2}. 
\end{proof}

\bigskip

\section{properties of monomial ideals  via polarization 
\hfill\break}\label{s:monomial-polarized} 

      For the purpose of all discussions in this section, unless
      otherwise stated, let $I$ be a monomial ideal in the polynomial
      ring $R=k[\xs]$ over a field $k$, whose polarization is the
      square-free monomial ideal $\PP(I)$ in the polynomial ring
      $$S=k[\x{i}{j} \ | \ 1 \leq i \leq n, 1 \leq j \leq a_i].$$ We assume
       that the polarizing sequence (as described in (\ref{e:sequence}) in
       Proposition~\ref{p:polarization}) is
       $$\vv=\nus$$  which is a regular
       sequence in $S/\PP(I)$ and
      $$R/I= S/\big (\PP(I) + (\vv) \big).$$

\subsection{Monomial ideals whose polarization is a simplicial tree}\label{ss:general-trees}

A natural question, and one that this paper is mainly concerned with,
is what properties of facet ideals of simplicial trees can be extended
to general (non-square-free) monomial ideals using polarization? In
other words, if for a monomial ideal $I$ in a polynomial ring $\PP(I)$
is the facet ideal of a tree (Definition~\ref{d:tree}), then what
properties of $\PP(I)$ are inherited by $I$?

The strongest tool when dealing with square-free monomial ideals is
induction-- either on the number of generators, or the number of variables in
the ambient polynomial ring. This is particularly the case when the facet
complex of the ideal is a tree, or in some cases when it just has a leaf. In
this section we show that via polarization, one can extend these tools to
monomial ideals in general. For a given monomial ideal $I$, we show that if
$\PP(I)$ is the facet ideal of a tree, and $\pp$ is a prime ideal containing
$I$, then $\PP(I_{\pp})$ is the facet ideal of a forest
(Theorem~\ref{t:localization}); this allows induction on number of
variables. Similarly, Theorem~\ref{t:joint-removal} provides us with a way to
use induction on number of generators of $I$.

\begin{theorem}[localization and polarization]\label{t:localization} 
If $\PP(I)$ is the facet ideal of a tree, and $\pp$ is a prime ideal
of $R$ containing $I$, then $\PP(I_{\pp})$ is the facet ideal of a
forest.
\end{theorem}

    \begin{proof} The first step is to show that it is enough to prove this
          for prime ideals of $R$ generated by a subset of $\{ \xs
          \}$. To see this, assume that $\pp$ is a prime ideal of $R$
          containing $I$, and that $\pp'$ is another prime of $R$
          generated by all $x_i \in \{ \xs \}$ such that $x_i \in \pp$
          (recall that the minimal primes of $I$ are generated by
          subsets of $\{ \xs \}$; see \cite[Corollary 5.1.5]{Vi2}). So
          $\pp' \subseteq \pp$. If $I=(\Ms)$, then
        $$I_{\pp'} = ({M_1}',\ldots,{M_q}')$$ where for each $i$, ${M_i}'$ is
        the image of $M_i$ in $I_{\pp'}$. In other words, ${M_i}'$ is obtained
        by dividing $M_i={x_1}^{a_1}\ldots {x_n}^{a_n}$ by the product of all
        the ${x_j}^{a_j}$ such that $x_j \notin \pp'$. But $x_j \notin \pp'$
        implies that $x_j \notin \pp$, and so it follows that ${M_i}' \in
        I_{\pp}$. Therefore $I_{\pp'} \subseteq I_{\pp}$. On the other hand since $\pp'
        \subseteq \pp$, $I_{\pp} \subseteq I_{\pp'}$, which implies that $I_{\pp'} =
        I_{\pp}$ (the equality and inclusions of the ideals here mean equality and
        inclusion of their generating sets).

         Now suppose $I=(\Ms)$, and $\pp=(x_1, \ldots, x_r)$ is a prime
         containing $I$. Suppose that for each $i$, we write
         $M_i=M_i'.M_i''$, where $$M_i' \in k[x_1, \ldots, x_r] {\rm \
         and \ } M_i'' \in k[x_{r+1}, \ldots, x_n]$$ so that
         $$I_{\pp}=(M_1',\ldots,M_t'),$$ where without loss of generality
         $M_1',\ldots,M_t'$ is a minimal generating set for $I_{\pp}$.
         
         We would like to show that the facet complex $\D$ of
         $\PP(I_{\pp})$ is a forest. Suppose that, again without loss of
         generality, $$I'=(\PP(M_1'), \ldots, \PP(M_s'))$$ is the facet
         ideal of a subcollection $\D'$ of $\D$. We need to show that
         $\D'$ has a leaf. 

         If $s=1$, then there is nothing to prove. Otherwise, suppose
         that $\PP(M_1)$ represents a leaf of the tree $\df(\PP(I))$,
         and $\PP(M_2)$ is a joint of $\PP(M_1)$. Then we have
         $$\PP(M_1) \cap \PP(M_i) \subseteq \PP(M_2) {\ for \ all \ }
         i\in \{2, \ldots, s\}.$$
         
        Now let $\x{e}{f}$ be in $\PP(M_1') \cap \PP(M_i')$ for some $i
        \in \{2, \ldots, s\}$. This implies that \begin{itemize}
        \item[($i$)] $\x{e}{f} \in \PP(M_1) \cap \PP(M_i) \subseteq
        \PP(M_2)$, and
        \item[($ii$)] $e \in \{1 \ldots, r\}$ \end{itemize} From ($i$)
        and ($ii$) we can conclude that $\x{e}{f} \in \PP(M_2')$,
        which proves that $\PP(M_1')$ is a leaf for $\D'$.  \end{proof}

\begin{remark}  It is not true in general that if $\pp$ is a (minimal)
 prime of $I$, then $\PP(I_{\pp})=\PP(I)_{\PP(\pp)}$. For example, if
$I=({x_1}^3,{x_1}^2x_2)$ and $\pp=(x_1)$, then $I_{\pp}=({x_1}^2)$ so
$\PP(I_{\pp})=(\x{1}{1}\x{1}{2})$, but $\PP(I)_{\PP(\pp)}=(\x{1}{1})$.
\end{remark}

       Another feature of simplicial trees is that they satisfy a
       generalization of K\"onig's theorem (\cite[Theorem~5.3]{F2}). 
       Below we explain how this property, and another
       property of trees that is very useful for induction, behave
       under polarization. 

       Recall that for a simplicial complex $\D$, $\al(\D)$ and
       $\beta(\D)$ are the vertex covering number and the independence
       number of $\D$, respectively (Definition~\ref{d:unmixed}). For
       simplicity of notation, if $I=(\Ms)$ is a monomial ideal, we
       let $\beta(I)$ denote the maximum cardinality of a subset of
       $\{\Ms\}$ consisting of pairwise coprime elements (so
       $\beta(\D)=\beta(\F(\D))$ for any simplicial complex $\D$).
 
       \begin{theorem}[joint removal and polarization]\label{t:joint-removal}
        Suppose $\Ms$ are monomials that form a minimal generating set
        for $I$, and $\PP(I)$ is the facet ideal of a simplicial
        complex $\D$. Assume that $\D$ has a leaf, whose joint
        corresponds to $\PP(M_1)$. Then, if we let $I'=(M_2, \ldots,
        M_q)$, we have $$\height I=\height I'.$$
       \end{theorem}

        \begin{proof} If $G$ is the joint of $\D$ corresponding to $\PP(M_1)$,
        then $\PP(I')=\F(\D\rmv{G})$. From \cite[Lemma~5.1]{F2} it
        follows that $\al(\D)=\al(\D \rmv{G})$, so that $\height \PP(I) =$
        $\height \PP(I')$, and therefore $\height I=\height I'$.  \end{proof}

       \begin{theorem}\label{t:konig} Suppose $\Ms$ are monomials that 
        form a minimal generating set for $I$, and $\PP(I)$ is the
        facet ideal of a simplicial tree $\D$. Then $\height
        I=\beta(I)$.
       \end{theorem}

       \begin{proof} We already know that $\height I=\height \PP(I) =
             \al(\D)$. It is also clear that $\be(I)=\be(\PP(I))$,
             since the monomials in a subset $\{M_{i_1}, \ldots,
             M_{i_r}\}$ of the generating set of $I$ are pairwise
             coprime if and only if the monomials in $\{\PP(M_{i_1}),
             \ldots, \PP(M_{i_r})\}$ are pairwise coprime. On the
             other hand, from \cite[Theorem~5.3]{F2} we know that
             $\al(\D)=\be(\D)$. Our claim follows immediately.  \end{proof}

We demonstrate how to apply these theorems via an example.

\begin{example} Suppose $I=({x_1}^3,{x_1}^2x_2x_3, 
{x_3}^2, {x_2}^3x_3)$.  Then
$$\PP(I)=(\x{1}{1}\x{1}{2}\x{1}{3},\x{1}{1}\x{1}{2}\x{2}{1}\x{3}{1},
\x{3}{1}\x{3}{2},\x{2}{1}\x{2}{2}\x{2}{3}\x{3}{1})$$ is the facet
ideal of the following simplicial complex (tree) $\D$.

\[ \includegraphics{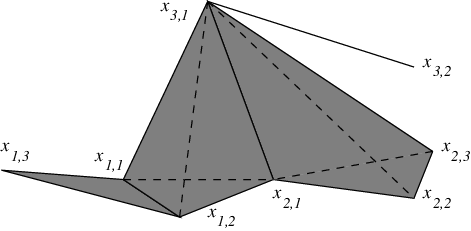} \]

Now $\al(\D)=\height I= 2$ because the prime of minimal height over
$I$ is $(x_1,x_3)$. From Theorem~\ref{t:konig} it follows that
$\be(I)=2$. This means that you can find a set of two monomials in the
generating set of $I$ that have no common variables: for example $\{
{x_1}^3, {x_3}^2\}$ is such a set.

Since the monomials ${x_1}^2x_2x_3$ and ${x_2}^3x_3$ polarize into
joints of $\D$, by Theorem~\ref{t:joint-removal} the ideals $$I, \  
({x_1}^3, {x_3}^2, {x_2}^3x_3), \
({x_1}^3,{x_1}^2x_2x_3, {x_3}^2), {\rm \ and \ }({x_1}^3, {x_3}^2)$$
all have the same height.

\end{example}

     We now focus on the Cohen-Macaulay property.  In \cite{F2} we
      showed that for a simplicial tree $\D$, $\F(\D)$ is a
      Cohen-Macaulay ideal if and only if $\D$ is an unmixed
      simplicial complex. The condition unmixed for $\D$ is equivalent
      to all minimal primes of the ideal $\F(\D)$ (which in this case
      are all the associated primes of $\F(\D)$) having the same
      height. In general, an ideal all whose associated primes have
      the same height (equal to the height of the ideal) is called an
      \emph{unmixed ideal}.

It now follows that

\begin{theorem}[Cohen-Macaulay criterion for trees] 
Let $\PP(I)$ be the facet ideal of a simplicial tree
$\D$. Then $R/I$ is Cohen-Macaulay if and only if $I$ is unmixed.
\end{theorem}

      \begin{proof} From Proposition~\ref{p:polarization}, $R/I$ is
            Cohen-Macaulay if and only if $S/\PP(I)$ is
            Cohen-Macaulay. By~\cite[Corollary~8.3]{F2}, this is
            happens if and only if $\PP(I)$ is
            unmixed. Corollary~\ref{c:polar-ass} now proves the claim.
            \end{proof}

	    If $\RR$ is a ring and $J$ is an ideal of $R$,
            then the Rees ring of $\RR$ along $J$ is defined as $$\RR[Jt]=
            \oplus_{n \in \mathbb N} {J^nt^n}.$$ Rees rings come up in
            the algebraic process of ``blowing up'' ideals. One reason
            that trees were defined as they are, is that their facet
            ideals produce normal and Cohen-Macaulay Rees rings
            (\cite{F1}).

 \begin{proposition}\label{p:Rees-CM} If $S[\PP(I)t]$ is
Cohen-Macaulay, then so is $R[It]$. Conversely, if we assume that $R$
and $S$ are localized at their irrelevant maximal ideals, then $R[It]$
being Cohen-Macaulay implies that $S[\PP(I)t]$ is Cohen-Macaulay.
\end{proposition}

      \begin{proof} Suppose that $\nus$ is the polarizing sequence as
      described before. For $i=1,\ldots,v-1$ let
      $$R_i=S/(\nu_1,\ldots,\nu_i), \ 
      I_i=\PP(I)/(\nu_1,\ldots,\nu_i), \ R_v=R {\rm \  and\ } I_v=I.$$

      Notice that $S[\PP(I)t]$ and $R[It]$ are both domains. Also note
      that for each $i$, $$S[\PP(I)t]/(\nu_1,\ldots,\nu_i)=R_i[I_it]$$ is
      the Rees ring of the monomial ideal $I_i$ in the polynomial ring
      $R_i$, and is therefore also a domain. Therefore $\nu_{i+1}$ is a
      regular element in the ring $S[\PP(I)t]/(\nu_1,\ldots,\nu_i)$, which
      means that $\nus$ is a regular sequence in $S[\PP(I)t]$.
       
      Similarly, we see that
      $$R[It]=S[\PP(I)t]/(\nus).$$  
 
      \cite[Theorem~2.1.3]{BH} now implies that if $S[\PP(I)t]$ is
      Cohen-Macaulay, then so is $R[It]$. The converse follows from
      \cite[Exercise~2.1.28]{BH}.\end{proof}

\begin{corollary}[Rees ring of a tree is Cohen-Macaulay]\label{c:Rees} 
Suppose that $\PP(I)$ is the facet ideal of a simplicial tree. Then the Rees
ring $R[It]$ of $I$ is Cohen-Macaulay.
\end{corollary} 
   
     \begin{proof} This follows from the Proposition~\ref{p:Rees-CM}, and
           from \cite[Corollary~4]{F1}, which states that the Rees ring
           of the facet ideal of a simplicial tree is
           Cohen-Macaulay. \end{proof}

\bigskip


\subsection{Polarization of sequentially Cohen-Macaulay ideals \hfill\break}

The main result of this section is that if the polarization of a
monomial ideal $I$ is the facet ideal of a tree, then $I$ is a
sequentially Cohen-Macaulay ideal. The theorem that implies this fact
(Proposition~\ref{p:duality}) is interesting in its own right. For a
square-free monomial ideal $J$, Eagon and Reiner \cite{ER} proved that
$J$ is Cohen-Macaulay if and only if its Alexander dual $\alx{J}$ has
a linear resolution. Herzog and Hibi~\cite{HH} then defined
componentwise linear ideals and generalized their result, so that a
square-free monomial ideal $J$ is sequentially Cohen-Macaulay if and
only if $\alx{J}$ is componentwise linear (see~\cite{HH}
or~\cite{F3}).  But even though Alexander duality has been generalized
to all monomial ideals from square-free ones, the criterion for
sequential Cohen-Macaulayness does not generalize: it is not true that
if $I$ is any monomial ideal, then $I$ is sequentially Cohen-Macaulay
if and only of $\alx{I}$ is a componentwise linear ideal; see
Miller~\cite{ezra-alex}. We show that the statement is true if
$\alx{I}$ is replaced by $\alx{\PP(I)}$.

\begin{definition}[\mbox{\cite[Chapter III, Definition~2.9]{S}}]\label{d:scm} 
 Let $M$ be a finitely generated ${\mathbb{Z}}$-graded module over a
  finitely generated ${\mathbb{N}}$-graded $k$-algebra, with $R_0=k$.
  We say that $M$ is \emph{sequentially Cohen-Macaulay} if there
  exists a finite filtration
  $$0=M_0 \subset M_1 \subset \ldots \subset M_r=M$$
  of $M$ by
  graded submodules $M_i$ satisfying the following two conditions.
\begin{enumerate}
\item[(a)] Each quotient $M_i/M_{i-1}$ is Cohen-Macaulay.

\item[(b)] $\dimn (M_1/M_0) < \dimn (M_2/M_1) < \ldots < \dimn (M_r/M_{r-1})$, 
where $\dimn$ denotes Krull dimension.
\end{enumerate} 
\end{definition}

We define a componentwise linear ideal in the square-free case using
\cite[Proposition~1.5]{HH}.

\begin{definition}[componentwise linear]\label{comp-def} 
  Let $I$ be a square-free monomial ideal in a polynomial ring
    $R$.  For a positive integer $k$, the $k$-th \emph{square-free
      homogeneous component} of $I$, denoted by $I_{[k]}$ is the ideal
    generated by all square-free monomials in $I$ of degree $k$.  The
    ideal $I$ above is said to be \emph{componentwise linear} if for
    all $k$, the square-free homogeneous component $I_{[k]}$ has a
    linear resolution.
\end{definition}

For our monomial ideal $I$, let $\QQ(I)$ denote the set of primary
ideals appearing in a reduced primary decomposition of $I$.  Suppose
that $h=\height I$ and $s= {\rm max}\{\height \qq \st \qq \in
\QQ(I)\}$, and set
 $$I_i = \bigcap_{\stackrel{\mbox{$\scriptstyle \qq \in \QQ(I)$}}{
\mbox{$\scriptstyle \height \qq \leq s-i$}}} \qq.$$

So we have the following filtration for $R/I$ (we assume that all
inclusions in the filtration are proper; if there is an equality
anywhere, we just drop all but one of the equal ideals).
\begin{eqnarray}\label{e:I-filtration}0=I=I_0 \subset I_1 \subset 
\ldots \subset I_{s-h} \subset R/I.\end{eqnarray}

If $R/I$ is sequentially Cohen-Macaulay, then by Theorem~\ref{t:main},
(\ref{e:I-filtration}) is the appropriate filtration that satisfies
Conditions (a) and (b) in Definition~\ref{d:scm}.

For the square-free monomial ideal $J=\PP(I)$, we similarly define
$$\QQ(J)=\{ {\rm minimal \ primes \ over \ } J \} {\rm \ and\ } J_i =
\bigcap_{\stackrel{\mbox{$\scriptstyle \pp \in \QQ(J)$}}{
\mbox{$\scriptstyle \height \pp \leq s-i$}}} \pp$$ where the numbers
$h$ and $s$ are the same as for $I$ because of
Proposition~\ref{p:polar-primary}.  It follows from
Proposition~\ref{p:polar-primary} and Corollary~\ref{c:polar-ass} that
for each $i$, $\PP(I_i)=J_i$ and the polarization sequence that
transforms $I_i$ into $J_i$ is a subsequence of $\vv=\nus$.

What we have done so far is to translate, via polarization, the
filtration (\ref{e:I-filtration}) of the quotient ring $R/I$ into one of
$S/J$:
\begin{eqnarray}\label{e:J-filtration} 0=J=J_0 \subset J_1 \subset 
\ldots \subset J_{s-h} \subset S/J.\end{eqnarray}

Now note that for a given $i$, the sequence $\vv$ is a
$J_{i+1}/J_i$-regular sequence in $S$, as $\vv$ is a regular sequence
in $S/J_i$, which contains $J_{i+1}/J_i$.  Also note that $$R/I_i
\simeq S/(J_i+\vv).$$  It follows that $J_{i+1}/J_i$ is
Cohen-Macaulay if and only if $I_{i+1}/I_i$ is (see
\cite[Exercise~2.1.27(c), Exercise~2.1.28, and Theorem~2.1.3]{BH}).

\begin{proposition}\label{p:duality} 
The monomial ideal $I$ is sequentially Cohen-Macaulay if and only if
 $\PP(I)$ is sequentially Cohen-Macaulay, or equivalently,
 $\alx{\PP(I)}$ is a componentwise linear ideal.
\end{proposition}

    \begin{proof} By \cite[Proposition~4.5]{F3}, $\alx{\PP(I)}$ is
           componentwise linear if and only if $\PP(I)$ is
           sequentially Cohen-Macaulay, which by the discussion above
           is equivalent to $I$ being sequentially Cohen-Macaulay.
    \end{proof}

\begin{corollary}[Trees are sequentially Cohen-Macaulay]\label{c:tree-scm}
Let $\PP(I)$ be the facet ideal of a simplicial tree. Then $I$ is
sequentially Cohen-Macaulay. \end{corollary}

    \begin{proof} This follows from Proposition~\ref{p:duality} and by
    \cite[Corollary~5.5]{F3}, which states that $\alx{\PP(I)}$ is a
    componentwise linear ideal.
    \end{proof} 

\bigskip 

\section{Further examples and remarks \hfill\break} 

To use the main results of this paper for computations on a given
monomial ideal, there are two steps. One is to compute the
polarization of the ideal, which as can be seen from the definition,
is a quick and simple procedure. This has already been implemented in
Macaulay2.  The second step is to determine whether the polarization
is the facet ideal of a tree, or has a leaf. Algorithms that serve
this purpose are under construction~\cite{CFS}.
 
\begin{remark}Let $I=(\Ms)$  be a monomial ideal in a polynomial 
ring $R$. If $\PP(I)$ is the facet ideal of a tree, then by
Corollary~\ref{c:Rees}, $R[It]$ is Cohen-Macaulay. But more is true:
if you drop any generator of $I$, for example if you consider
$I'=(M_1,\ldots,\hat{M_i}, \ldots, M_q)$, then $R[I't]$ is still
Cohen-Macaulay. This is because $\PP(I')$ corresponds to the facet
ideal of a forest, so one can apply the same result.  
\end{remark}

A natural question is whether one can say the same with the property
``Cohen-Macaulay'' replaced by ``normal''. If $I$ is square-free, this
is indeed the case. But in general, polarization does not preserve
normality of ideals.

\begin{example}[normality and polarization]\label{r:normal-polar}
A valid question is whether Proposition~\ref{p:Rees-CM} holds if
the word ``Cohen-Macaulay'' is replaced with ``normal'', given that
simplicial trees have normal facet ideals (\cite{F1})?

The answer is negative. Here is an example. 
 
Let $I=({x_1}^3, {x_1}^2{x_2},{x_2}^3)$ be an ideal of
$k[x_1,x_2]$. Then $I$ is not normal; this is because $I$ is not even
integrally closed: ${x_1}{x_2}^2 \in \overline{I}$ as
$({x_1}{x_2}^2)^3 -{x_1}^3{x_2}^6=0$, but ${x_1}{x_2}^2 \not \in I$.
Now $$\PP(I)=(\x{1}{1}\x{1}{2}\x{1}{3},
\x{1}{1}\x{1}{2}\x{2}{1},\x{2}{1}\x{2}{2}\x{2}{3})$$ is the facet
ideal of the tree

\[ \includegraphics{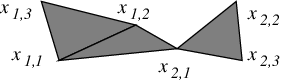} \]

 which is normal by~\cite{F1}.
\end{example}

The reason that normality (or integral closure in general) does not
pass through polarization is much more basic: polarization does not
respect multiplication of ideals, or monomials. Take, for example, two
monomials $M$ and $N$ and two monomial ideals $I$ and $J$, such that
$MN \in IJ$. It is not necessarily true that $\PP(M)\PP(N) \in
\PP(I)\PP(J)$.

Indeed, let $I=J=(x_1x_2)$ and $M={x_1}^2$ and $N={x_2}^2$. Then
$MN={x_1}^2{x_2}^2 \in IJ$. But $$\PP(M)=\x{1}{1}\x{1}{2}, \ 
\PP(N)=\x{2}{1}\x{2}{2}, {\rm \ and \ } \PP(I)=\PP(J)=(\x{1}{1}\x{2}{1})$$ and
clearly $\PP(M)\PP(N) \notin \PP(I)\PP(J)=({\x{1}{1}}^2{\x{2}{1}}^2)$.

\begin{remark}It is useful to think of polarization as a chain of 
substitutions. This way, as a monomial ideal $I$ gets polarized, the
ambient ring extends one variable at a time. All the in-between ideals
before we hit the final square-free ideal $\PP(I)$ have the same
polarization. 

      For example let $J=(x^2,xy,y^3) \subseteq k[x,y]$.  We use a
      diagram to demonstrate the process described in the previous
      paragraph.  Each linear form $a-b$ stands for ``replacing the
      variable $b$ with $a$'', or vice versa, depending on which
      direction we are going.

      $$ J=(x^2,xy,y^3) \stackrel{u-x}{\longrightarrow} J_1=(xu,uy,y^3)
         \stackrel{v-y}{\longrightarrow} J_2=(xu,uv,y^2v)$$
      $$ \stackrel{w-y}{\longrightarrow} J_3=(xu,uv,yvw)\ .$$

      This final square-free monomial ideal $J_3$ is the polarization
      of $J$, and is isomorphic to $\PP(J)$ as we defined it in
      Definition~\ref{d:polarization}. Note that all ideals $J$,
      $J_1$, $J_2$ and $J_3$ have the same (isomorphic) polarization.
      We can classify monomial ideals according to their
      polarizations. An interesting question is to see what properties
      do ideals in the same polarization class have. A more difficult
      question is how far can one ``depolarize'' a square-free
      monomial ideal $I$, where by depolarizing $I$ we mean finding
      monomial ideals whose polarization is equal to $I$, or
      equivalently, traveling the opposite direction on the above
      diagram.

     \end{remark}

\appendix

\section{Appendix: Primary decomposition in a sequentially \\
  Cohen-Macaulay module}\label{appendix}

  The purpose of this appendix is to study, using basic facts about
  primary decomposition of modules, the structure of the submodules
  appearing in the (unique) filtration of a sequentially
  Cohen-Macaulay module $M$. The main result (Theorem~\ref{t:main})
  states that each submodule appearing in the filtration of $M$ is the
  intersection of all primary submodules whose associated primes have
  a certain height and appear in an irredundant primary decomposition
  of the $0$-submodule of $M$. Similar results, stated in a different
  language, appear in~\cite{SC}; the author thanks J\"urgen Herzog for
  pointing this out.

We first record two basic lemmas that we shall use later (the second
one is an exercise in Bourbaki~\cite{B}). Throughout the discussions
below, we assume that $R$ is a finitely generated algebra over a
field, and $M$ is a finite module over $R$.

\begin{lemma}\label{l:decomposition} Let $\QQs{1}, \ldots, \QQs{t}, \PP$ 
all be primary submodules of an $R$-module $M$, such that
$\ass{M/\QQs{i}}=\{\qq_i\}$ and $\ass{M/\PP}=\{\wp\}$. If $\QQs{1} \cap \ldots
\cap \QQs{t} \subseteq \PP$ and $\QQs{i} \not\subseteq \PP$ for some $i$, then  
there is a $j \neq i$ such that $\qq_j \subseteq \wp$.
\end{lemma}

  \begin{proof} Let $x \in \QQs{i} \setminus \PP$. For each $j$ not equal to
  $i$, pick the positive integer $m_j$ such that $\qq_j^{m_j} x
  \subseteq \QQs{j}.$ So we have that $$\qq_1^{m_1}\ldots
  \qq_{i-1}^{m_{i-1}} \qq_{i+1}^{m_{i+1}}\ldots \qq_t^{m_t}x \subseteq
  \QQs{1} \cap \ldots \cap \QQs{t} \subseteq \PP \ \Longrightarrow \
  \qq_1^{m_1}\ldots \qq_{i-1}^{m_{i-1}} \qq_{i+1}^{m_{i+1}}\ldots
  \qq_t^{m_t} \subseteq \wp$$ where the second inclusion is because $x
  \notin \PP$. Hence for some $j\neq i$, $\qq_j \subseteq \wp$.  \end{proof}

\begin{lemma}\label{l:exercise} Let $M$ be an $R$-module and
 $N$ be a submodule of $M$. Then for every $\wp \in \ass{M/N}$, if $\wp
 \not\supseteq \ann{N}$, then $\wp \in \ass{M}$.
\end{lemma}

    \begin{proof} Since $\wp \in \ass{M/N}$, there exists $x \in M \setminus
     N$ such that $\wp = \ann{x}$; in other words $\wp x \subseteq N.$
     Suppose $\ann{N} \not\subseteq \wp$, and let $y \in \ann{N}
     \setminus \wp$. Now $y\wp x =0$, and so $\wp \subseteq \ann{yx}$
     in $M$.  On the other hand, if $z \in \ann{yx}$, then $zyx=0
     \subseteq N$ and so $zy \in \wp$. But $y \notin \wp$, so $z \in
     \wp$.  Therefore $\wp \in \ass{M}$.

    \end{proof} 

Suppose $M$ is a \scm module with filtration as in Definition~\ref{d:scm}.
We adopt the following notation. For a given integer $j$, we let
$$\ass{M}_j=\{\wp \in \ass{M} \st \height \wp =j\}.$$ Suppose that all
the $j$ where $\ass{M}_j \neq \emptyset$ form the sequence of
integers $$0 \leq h_1 < \ldots < h_c \leq \dimn R$$ so that $\ass{M}=
\bigcup_{1 \leq j \leq c} \ass{M}_{h_j}.$

\begin{proposition}\label{p:work} For all $i=0,\ldots,r-1$, we have 
\begin{enumerate}

\item\label{ii:0} $\ass{M_{i+1}/M_i} \cap \ass{M} \neq \emptyset$;

\item\label{ii:1} $\ass{M}_{h_{r-i}} \subseteq \ass{M_{i+1}/M_i}$ and
$c=r$;

\item\label{ii:2} If $\wp \in \ass{M_{i+1}}$, then $\height \wp \geq
h_{r -i}$;

\item\label{ii:3} If $\wp \in \ass{M_{i+1}/M_i}$, then $\ann{M_i}
\not\subseteq \wp$;

\item\label{ii:4} $\ass{M_{i+1}/M_i} \subseteq \ass{M}$;

\item\label{ii:5} $\ass{M_{i+1}/M_i}= \ass{M}_{h_{r-i}}$;

\item\label{ii:6} $\ass{M/M_i} = \ass{M}_{\leq h_{r-i}}$;

\item\label{ii:7} $\ass{M_{i+1}} = \ass{M}_{\geq h_{r-i}}$.

\end{enumerate}
\end{proposition}

\begin{proof}
\begin{enumerate}

  \item We use induction on the length $r$ of the filtration of $M$.
            The case $r=1$ is clear, as we have a filtration $0
            \subset M$, and the assertion follows. Now suppose the
            statement holds for sequentially Cohen-Macaulay modules
            with filtrations of length less than $r$. Notice that
            $M_{r-1}$ that appears in the filtration of $M$ in
            Definition~\ref{d:scm} is also sequentially
            Cohen-Macaulay, and so by the induction hypothesis, we
            have $$\ass{M_{i+1}/M_i} \cap \ass{M_{r-1}} \neq \emptyset
            {\rm \ for \ } i=0,\ldots,r-2$$ and since $\ass{M_{r-1}}
            \subseteq \ass{M}$ it follows that $$\ass{M_{i+1}/M_i}
            \cap \ass{M} \neq \emptyset {\rm \ for \ }
            i=0,\ldots,r-2.$$

            It remains to show that $\ass{M/M_{r-1}}\cap \ass{M} \neq
            \emptyset$.

           For each $i$, $M_{i-1} \subset M_i$, so
            we have (\cite{B} Chapter~IV) 
                        \begin{eqnarray}\label{e:first}
            \ass{M_1} \subseteq \ass{M_2} \subseteq \ass{M_1} \cup
            \ass{M_2/M_1}\end{eqnarray}

           The inclusion $M_2 \subseteq M_3$ along with the inclusions
           in (\ref{e:first}) imply that
            $$\begin{array}{ll}
             \ass{M_2} \subseteq \ass{M_3}\hspace{-.1in}& 
              \subseteq \ass{M_2} \cup \ass{M_3/M_2} \\
           & \subseteq \ass{M_1} \cup \ass{M_2/M_1}\cup \ass{M_3/M_2}. 
	    \end{array}$$
          
          If we continue this process inductively, at the $i$-th stage
          we have
          $$\begin{array}{ll}
          \ass{M_i}\hspace{-.1in}&\subseteq \ass{M_{i-1}} \cup 
          \ass{M_i/M_{i-1}}\\
          &\subseteq \ass{M_1} \cup \ass{M_2/M_1}\cup \ass{M_3/M_2}
             \cup \ldots \cup \ass{M_i/M_{i-1}}\end{array}$$ 
          and finally, when  $i=r$ it gives
          \begin{eqnarray}\label{e:final} \ass{M}
            \subseteq \ass{M_1} \cup \ass{M_2/M_1}
             \cup \ldots \cup \ass{M/M_{r-1}}.\end{eqnarray}
 
         Because of Condition (b) in Definition~\ref{d:scm}, and the
         fact that each $M_{i+1}/M_i$ is Cohen-Macaulay (and hence all
         its associated primes have the same height; see~\cite{BH}
         Chapter~2), if for every $i$ we pick $\wp_i \in
         \ass{M_{i+1}/M_i}$, then $$ h_c \geq \height \wp_0 > \height
         \wp_1 > \ldots > \height \wp_{r-1}.$$ where the
         left-hand-side inequality comes from the fact that $\ass{M_1}
         \subseteq \ass{M}$. By our induction hypothesis, $\ass{M}$
         intersects $\ass{M_{i+1}/M_i}$ for all $i \leq r-2$, and so
         because of (\ref{e:final}) we conclude that $$\height \wp_i =
         h_{c-i}, {\rm \ and \ } \ass{M}_{h_{c-i}} \subseteq
         \ass{M_{i+1}/M_i} {\rm \ for \ } 0 \leq i \leq r-2.$$ And now
         $\ass{M}_{h_0}$ has no choice but to be included in
         $\ass{M/M_{r-1}}$, which settles our claim. It also follows
         that $c=r$.

 \item See the proof for part~\ref{ii:0}.

 \item We use induction. The case $i=0$ is clear, since for
       every $\wp \in \ass{M_1}=\ass{M_1/M_0}$ we know
       from part~\ref{ii:1} that $\height \wp = h_r$.  Suppose the
       statement holds for all indices up to $i-1$. Consider the
       inclusion $$\ass{M_{i}} \subseteq \ass{M_{i+1}} \subseteq
       \ass{M_i} \cup \ass{M_{i+1}/M_i}.$$ From part~\ref{ii:1} and the
       induction hypothesis it follows that if $\wp \in \ass{M_{i+1}}$
       then $\height \wp \geq h_{r-i}$. 

\item Suppose $\ann{M_i} \subseteq \wp$. Since $\sqrt{\ann{M_i}}=
       \bigcap_{\wp' \in \ass{M_i}} \wp'$, it follows that
       $\bigcap_{\wp' \in \ass{M_i}} \wp' \subseteq \wp,$ so there is
       a $\wp' \in \ass{M_i}$ such that $\wp' \subseteq \wp$. But from
       parts~\ref{ii:1} and~\ref{ii:2} above it follows that $\height
       \wp' \geq h_{r-i+1}$ and $\height \wp =h_{r-i}$, which is a
       contradiction.

\item From part~\ref{ii:3} and Lemma~\ref{l:exercise}, it follows that
        $\ass{M_{i+1}/M_i} \subseteq \ass{M_{i+1}} \subseteq \ass{M}.$

\item  This follows from parts~\ref{ii:1} and~\ref{ii:4}, and  the fact that
       $M_{i+1}/M_i$ is Cohen-Macaulay, and hence all associated
       primes have the same height.

\item We show this by induction on $e=r-i$. The case $e=1$ (or
      $i=r-1$) is clear, because by part~\ref{ii:5} we have
      $\ass{M/M_{r-1}}=\ass{M}_{h_1}=\ass{M}_{\leq h_1}.$

      Now suppose the equation holds for all integers up to $e-1$
      (namely $i=r-e+1$), and we would like to prove the statement for
      $e$ (or $i=r-e$).  Since $M_{i+1}/M_i \subseteq M/M_i$, we have
      \begin{eqnarray}\label{e:third} \ass{M_{i+1}/M_i} \subseteq 
      \ass{M/M_i} \subseteq \ass{M_{i+1}/M_i} \cup
      \ass{M/M_{i+1}}\end{eqnarray}

     By the induction hypothesis and part~\ref{ii:5} we know that
      $$\ass{M/M_{i+1}}=\ass{M}_{\leq h_{r-i-1}} {\rm \ and \ }
      \ass{M_{i+1}/M_i}= \ass{M}_{h_{r-i}},$$ which put together with
      (\ref{e:third}) implies that
      $$\ass{M}_{h_{r-i}} \subseteq \ass{M/M_i} \subseteq
     \ass{M}_{\leq h_{r-i}}$$ We still have to show that $\ass{M/M_i}
     \supseteq \ass{M}_{\leq h_{r-i-1}}$.

     Let $$\wp \in \ass{M}_{\leq {h_{r-i-1}}}= \ass{M/M_{i+1}}=
     \ass{(M/M_{i})/(M_{i+1}/M_i)}.$$ If $\wp \supseteq
     \ann{M_{i+1}/M_i}$, then (by part~\ref{ii:5}) $$\wp \supseteq
     \bigcap_{\qq \in \ass{M}_{h_{r-i}}} \qq \ \Longrightarrow\ \wp
     \supseteq \qq {\rm \ for\ some\ } \qq \in \ass{M}_{h_{r-i}}$$ which
     is a contradiction, as $\height \wp \leq h_{r-i-1} < \height \qq$.
     It follows from Lemma~\ref{l:exercise} that $\wp \in
     \ass{M/M_i}$.

\item The argument is based on induction, and exactly the same as the
      one in part~\ref{ii:3}, using more information; from the
      inclusions
      $$\ass{M_{i}} \subseteq \ass{M_{i+1}} \subseteq \ass{M_i} \cup
      \ass{M_{i+1}/M_i},$$ the induction hypothesis, and
      part~\ref{ii:5} we deduce that
      $$ \ass{M}_{\geq h_{r-i+1}} \subseteq \ass{M_{i+1}} \subseteq
      \ass{M}_{\geq h_{r-i+1}} \cup \ass{M}_{h_{r-i}}, $$ which put
      together with part~\ref{ii:3}, along with Lemma~\ref{l:exercise}
      produces the equality.

\end{enumerate} 
\end{proof}

Now suppose that as a submodule of $M$, $M_0=0$ has an irredundant
primary decomposition of the form:
\begin{eqnarray}\label{e:decomposition} 
M_0=0=\bigcap_{1 \leq j \leq r} \QQs{1}^{h_j} \cap \ldots \cap
\QQs{s_j}^{h_j}\end{eqnarray}
where for a fixed $j \leq r$ and $e
\leq s_j$, $\QQs{e}^{h_j}$  is a  primary submodule of $M$
with $\ass{M/\QQs{e}^{h_j}}=\{\wp_e^{h_j}\}$ and
$\ass{M}_{h_j}=\{\wp_1^{h_j}, \ldots, \wp_{s_j}^{h_j}\}.$

\begin{theorem}\label{t:main} Let $M$ be a sequentially Cohen-Macaulay 
module with filtration as in Definition~\ref{d:scm}, and suppose that
$M_0=0$ has a primary decomposition as in (\ref{e:decomposition}).
Then for each $i=0,\ldots,r-1$, $M_i$ has the following primary
decomposition
\begin{eqnarray}\label{e:i-decomposition}
M_i=\bigcap_{1 \leq j \leq r-i} \QQs{1}^{h_j} \cap \ldots \cap
\QQs{s_j}^{h_j}.\end{eqnarray}
\end{theorem}

     \begin{proof} We prove this by induction on $r$ (length of the
       filtration). The case $r=1$ is clear, as the filtration is of
       the form $0=M_0 \subset M$. Now consider $M$ with filtration 
       $$0=M_0 \subset M_1 \subset \ldots \subset M_r=M.$$ Since
       $M_{r-1}$ is a sequentially Cohen-Macaulay module of length
       $r-1$, it satisfies the statement of the theorem.  We first
       show that $M_{r-1}$ has a primary decomposition as described in
       (\ref{e:i-decomposition}). From part~\ref{ii:6} of
       Proposition~\ref{p:work} it follows that
       $$\ass{M/M_{r-1}}=\ass{M}_{h_1}$$ and so for some
       $\wp_e^{h_1}$-primary submodules $\PP_e^{h_1}$ of $M$ ($1 \leq
       e \leq s_j$), we have
       \begin{eqnarray}\label{e:induction-decomposition}M_{r-1}=
       \PP_1^{h_1} \cap \ldots \cap \PP_{s_1}^{h_1}.\end{eqnarray}
       
      We would like to show that $\QQs{e}^{h_1}=\PP_e^{h_1}$ for
      $e=1,\ldots,s_1$.

       Fix $e=1$ and assume $\QQs{1}^{h_1} \not \subset \PP_1^{h_1}$. From the
       inclusion $M_0 \subset \PP_1^{h_1}$ and
       Lemma~\ref{l:decomposition} it follows that for some $e$ and
       $j$ (with $e\neq 1$ if $j=1$), we have $\wp_e^{h_j} \subseteq
       \wp_1^{h_1}$. Because of the difference in heights of these
       ideals the only conclusion is $\wp_e^{h_j}=\wp_1^{h_1}$, which
       is not possible. With a similar argument we deduce that
       $\QQs{e}^{h_1} \subset \PP_e^{h_1}$, for $e=1, \ldots, s_1$.

       Now fix $j \in \{1,\ldots,r\}$ and $e \in \{1,\ldots,
       s_j\}$. If $M_{r-1}=\QQs{e}^{h_j}$ we are done.  Otherwise, note
       that for every $j$ and $\wp_e^{h_j}$-primary submodule
       $\QQs{e}^{h_j}$ of $M$,
       $$\QQs{e}^{h_j} \cap M_{r-1}$$ is a $\wp_e^{h_j}$-primary
       submodule of $M_{r-1}$ (as $\emptyset \neq
       \ass{M_{r-1}/(\QQs{e}^{h_j}\cap M_{r-1})} =
       \ass{(M_{r-1}+ \QQs{e}^{h_j})/\QQs{e}^{h_j}} \subseteq
       \ass{M/\QQs{e}^{h_j}}=\{\wp_e^{h_j}\}$). So $M_0=0$ as a
       submodule of $M_{r-1}$ has a primary decomposition
      $$M_0 \cap M_{r-1}=0= \bigcap_{1 \leq j \leq r} (\QQs{1}^{h_j}
       \cap M_{r-1}) \cap \ldots \cap (\QQs{s_j}^{h_j} \cap
       M_{r-1}).$$ From Proposition~\ref{p:work} part~\ref{ii:7} it
       follows that $$ \ass{M_{r-1}} = \ass{M}_{\geq h_2}$$ so the
       components $\QQs{t}^{h_1} \cap M_{r-1}$ are redundant for $t=1,
       \ldots, s_1$, so for each such $t$ we have
       $$\bigcap_{\QQs{e}^{h_j} \neq \QQs{t}^{h_1}} (\QQs{1}^{h_j} \cap
       M_{r-1}) \subseteq \QQs{t}^{h_1} \cap M_{r-1}.$$ If
       $\QQs{e}^{h_j} \cap M_{r-1} \not \subseteq \QQs{t}^{h_1} \cap
       M_{r-1}$ for some $e$ and $j$ (with $\QQs{e}^{h_j} \neq
       \QQs{t}^{h_1}$), then by Lemma~\ref{l:decomposition} for some
       such $e$ and $j$ we have $\wp_e^{h_j} \subseteq \wp_t^{h_1}$,
       which is a contradiction (because of the difference of
       heights).

       Therefore, for each $t$ ($1 \leq t \leq s_1$), there exists
       indices $e$ and $j$ (with $\QQs{e}^{h_j} \neq \QQs{t}^{h_1}$) such
       that $$\QQs{e}^{h_j} \cap M_{r-1} \subseteq \QQs{t}^{h_1} \cap
       M_{r-1}.$$ It follows now, from the primary decomposition of
       $M_{r-1}$ in (\ref{e:induction-decomposition}) that for a fixed $t$
       $$\PP_1^{h_1} \cap \ldots \cap \PP_{s_1}^{h_1} \cap \QQs{e}^{h_j}
       \subseteq \QQs{t}^{h_1}.$$

     Assume $\PP_t^{h_1} \not \subseteq \QQs{t}^{h_1}$. Applying
     Lemma~\ref{l:decomposition} again, we deduce that $$\wp_{e}^{h_j}
     \subseteq \wp_t^{h_1}, {\rm \ or\ there\ is\ } t'\neq t {\rm \
     such\ that\ } \wp_{t'}^{h_1} \subseteq \wp_t^{h_1}.$$ Neither of
     these is possible, so $\PP_t^{h_1} \subseteq \QQs{t}^{h_1}$ for all
     $t$.

      We have therefore proved that $$M_{r-1}=\QQs{1}^{h_1} \cap \ldots \cap
       \QQs{s_1}^{h_1}.$$

       By the induction hypothesis, for each $i \leq r-2$, $M_i$ has
       the following primary decomposition
      $$ M_i = \bigcap_{2 \leq j \leq r-i} (\QQs{1}^{h_j} \cap M_{r-1}) \cap
                \ldots \cap (\QQs{s_j}^{h_j} \cap M_{r-1}) = \bigcap_{1 \leq j
                \leq r-i} \QQs{1}^{h_j} \cap \ldots \cap \QQs{s_j}^{h_j}$$ 
      which proves the theorem.

     \end{proof}


\end{document}